\documentclass[12pt]{article}
\usepackage{amsmath}
\usepackage{amssymb}
\usepackage{latexsym}
\usepackage{amsthm}
\usepackage{mathrsfs}
\usepackage{enumitem}
\usepackage[colorlinks,
            linkcolor=red,
            anchorcolor=blue,
            citecolor=green]{hyperref}

\newtheorem{thm}{Theorem}[section]

\newtheorem{prop}[thm]{Proposition}

\newtheorem{conj}[thm]{Conjecture}

\begin{document}
\begin{center}
{\large \bf  Log-concavity of the Fennessey-Larcombe-French Sequence}
\end{center}
\begin{center}
Arthur L.B. Yang$^{1}$ and James J.Y. Zhao$^{2}$\\[6pt]

$^{1}$Center for Combinatorics, LPMC-TJKLC\\
Nankai University, Tianjin 300071, P. R. China\\[6pt]

$^{2}$Center for Applied Mathematics\\
Tianjin University, Tianjin 300072, P. R. China\\[8pt]

Email: $^{1}${\tt yang@nankai.edu.cn},
       $^{2}${\tt jjyzhao@tju.edu.cn}
\end{center}

\noindent\textbf{Abstract.}
We prove the log-concavity of the Fennessey-Larcombe-French sequence based on its three-term recurrence relation, which was recently conjectured by Zhao. The key ingredient of our approach is
a sufficient condition for log-concavity of a sequence subject to certain three-term recurrence.

\noindent \emph{AMS Classification 2010:} Primary 05A20


\section{Introduction}

The objective of this paper is to prove the log-concavity conjecture of the Fennessey-Larcombe-French sequence, which was posed by Zhao \cite{zhaocf} in the study of log-balancedness of combinatorial sequences.

Let us begin with an overview of Zhao's conjecture.
Recall that a sequence $\{a_k\}_{k \geq 0}$ is said to be {log-concave} if
$$a_{k}^2\geq a_{k+1}a_{k-1}, \quad  \mbox{for $k\geq 1$,} $$
and it is log-convex if
$$a_{k}^2\leq a_{k+1}a_{k-1}, \quad  \mbox{for $k\geq 1$.} $$
We say that $\{a_k\}_{k \geq 0}$ is log-balanced if the sequence itself is log-convex while $\{\frac{a_k}{k!}\}_{k \geq 0}$ is log-concave.

The Fennessey-Larcombe-French sequence $\{V_n\}_{n\geq 0}$ can be given by the following three-term recurrence relation \cite{lff2003}
\begin{align}\label{eq-rec}
n(n+1)^2 V_{n+1}=8n(3n^2+5n+1)V_n-128(n-1)(n+1)^2V_{n-1}, \quad \mbox{for }n\geq 1,
\end{align}
with the initial values $V_0=1$ and $V_1=8$.
This sequence was introduced by Larcombe, French and Fennessey \cite{lff2002}, in connection with
a series expansion of the complete elliptic integral of the second kind, precisely,
\begin{align*}
\int_0^{\pi/2}\sqrt{1-c^2\sin^2\theta}\,d\theta
=\frac{\pi\sqrt{1-c^2}}{2}\sum\limits_{n=0}^\infty \left(\frac{1-\sqrt{1-c^2}}{16}\right)^n V_n.
\end{align*}

The Fennessey-Larcombe-French sequence is closely related to the Catalan-Larcombe-French sequence, which was first studied by E. Catalan \cite{Catalan} and later examined and clarified by Larcombe and French \cite{lf2000}.
Let $\{P_n\}_{n\geq 0}$ denote the Catalan-Larcombe-French sequence, and the following three-term recurrence relation holds:
 \begin{align*}
 (n+1)^2P_{n+1}=8(3n^2+3n+1)P_{n}-128n^2P_{n-1}, \quad\mbox{for } n\geq 1,
 \end{align*}
 with $P_0=1$ and $P_1=8$.
 As a counterpart of $V_n$, the numbers $P_n$ appear as coefficients in the series expansion of the complete elliptic integral of the first kind, precisely,
\begin{align*}
\int_0^{\pi/2}\frac{1}{\sqrt{1-c^2\sin^2\theta}}\,d\theta
=\frac{\pi}{2}\sum\limits_{n=0}^\infty \left(\frac{1-\sqrt{1-c^2}}{16}\right)^n P_n.
\end{align*}
Many interesting properties have been found for the Catalan-Larcombe-French sequence and the Fennessey-Larcombe-French sequence, and the reader may consult references \cite{jlf2004, jv2010, lf2000, lff2002, lff2003, Zagier}.

Recently, there has arisen an interest in the study of the log-behavior of the Catalan-Larcombe-French sequence. For instance, Xia and Yao \cite{XY2013} obtained the log-convexity of the Catalan-Larcombe-French sequence, and confirmed a conjecture of Sun \cite{sun2013}. By using a log-balancedness criterion due to Do\v{s}li\'{c} \cite{Doslic}, Zhao \cite{zhao2014} proved the log-balancedness of the Catalan-Larcombe-French sequence.

Zhao further studied the log-behavior of the Fennessey-Larcombe-French sequence, and obtained the following result.
\begin{thm}[\cite{zhaocf}]\label{thm-z}
Both $\{nV_n\}_{n\geq 1}$ and $\{\frac{V_n}{(n-1)!}\}_{n\geq 1}$ are log-concave.
\end{thm}
She also made the following conjecture \cite{zhaocf,zhaopr}.
\begin{conj}\label{conj-z}
The Fennessey-Larcombe-French sequence $\{V_n\}_{n\geq 1}$ is log-concave.
\end{conj}
Note that the Hadamard product of two log-concave sequences without internal zeros is still log-concave, see \cite[Proposition 2]{stanley1989}. Since both $\{n\}_{n\geq 1}$ and $\{\frac{1}{(n-1)!}\}_{n\geq 1}$ are log-concave, Conjecture \ref{conj-z} implies Theorem \ref{thm-z}.

In this paper, we obtain a sufficient condition for proving the log-concavity of a sequence satisfying a three-term recurrence. Then we give an affirmative answer to Conjecture \ref{conj-z} by using this criterion. By further employing a result of Wang and Zhu \cite[Theorem 2.1]{WangZhu}, we derive the monotonicity of the sequence $\{\sqrt[n]{V_{n+1}}\}_{n\geq 1}$ from the log-concavity of $\{{V_{n}}\}_{n\geq 1}$.

\section{Log-concavity derived from three-term recurrence}

The aim of this section is to prove the log-concavity of the Fennessey-Larcombe-French sequence based on its three-term recurrence relation.

We first give a sufficient condition for log-concavity of a positive sequence subject to certain three-term recurrence. It should be mentioned that the log-behavior of sequences satisfying three-term recurrences has been extensively studied, see Liu and Wang \cite{LiuWang}, Chen and Xia \cite{ChenXia}, Chen, Guo and Wang \cite{CGW2014}, and Wang and Zhu \cite{WangZhu}. However, most of these studies have focused on the log-convexity of such sequences instead of their log-concavity.
Our criterion for determining the log-concavity of a sequence satisfying a three-term recurrence is as follows.

\begin{prop}\label{criterion}
Let $\{S_n\}_{n\geq 0}$ be a positive sequence satisfying the following recurrence relation:
\begin{align}\label{eq-S}
a(n)S_{n+1}+b(n)S_n+c(n)S_{n-1}=0, \quad \mbox{for $n\geq 1$,}
\end{align}
where $a(n),b(n)$ and $c(n)$ are real functions in $n$. Suppose that there exists an integer $n_0$ such that for any $n> n_0$,
\begin{itemize}
\item[$(i)$] it holds $a(n)>0$, and
\item[$(ii)$] either $b^2(n)<4a(n)c(n)$ or $\frac{S_n}{S_{n-1}}\geq \frac{-b(n)+\sqrt{b^2(n)-4a(n)c(n)}}{2a(n)}$.
\end{itemize}
Then the sequence $\{S_n\}_{n\geq n_0}$ is log-concave, namely, $S_n^2\geq S_{n+1}S_{n-1}$ for any $n>n_0$.
\end{prop}

\proof Let $r(n)=\frac{S_n}{S_{n-1}}$. It suffices to show that $r(n)\geq r(n+1)$ for any $n>n_0$.
On one hand, the conditions $(i)$ and $(ii)$ imply that
\begin{align*}
a(n)r^2(n)+b(n)r(n)+c(n)\geq 0, \quad \mbox{for $n> n_0$.}
\end{align*}
Since $\{S_n\}_{n\geq 0}$ is a positive sequence, so is $\{r_n\}_{n\geq 1}$.
Thus, the above inequality is equivalent to the following
\begin{align}\label{eq-1}
a(n)r(n)+b(n)+\frac{c(n)}{r(n)}\geq 0, \quad \mbox{for $n> n_0$.}
\end{align}
On the other hand, dividing both sides of \eqref{eq-S} by $S_n$, we obtain
\begin{align}\label{eq-2}
a(n)r(n+1)+b(n)+\frac{c(n)}{r(n)}=0.
\end{align}
Combining \eqref{eq-1} and \eqref{eq-2}, we get
\begin{align*}
a(n)r(n+1)\leq a(n)r(n), \quad \mbox{for $n> n_0$.}
\end{align*}
By the condition $(i)$, we have $r(n+1)\leq r(n)$ for any $n> n_0$. This completes the proof.
\qed

We are now able to give the main result of this section, which offers an affirmative answer to Conjecture \ref{conj-z}.

\begin{thm}\label{thm-main} Let $\{V_n\}_{n\geq 0}$ be the Fennessey-Larcombe-French sequence given by \eqref{eq-rec}. Then, for any $n\geq 2$, we have $V_n^2\geq V_{n-1}V_{n+1}$.
\end{thm}

\proof By the recurrence relation \eqref{eq-rec}, we have
$V_1=8, V_2=144, V_3=2432$ and $V_4=40000$. It is easy to verify that
$V_2^2\geq V_1V_3$ and $V_3^2\geq V_2V_4$.

We proceed to use Proposition \ref{criterion}
to prove that $V_n^2> V_{n-1}V_{n+1}$ for $n>3$, namely taking $n_0=3$.
For the sequence $\{V_n\}_{n\geq 0}$, the corresponding polynomials $a(n),b(n),c(n)$ appearing in Proposition \ref{criterion} are as follows:
\begin{align*}
a(n)&=n(n+1)^2,\\[5pt]
b(n)&=-8n(3n^2+5n+1),\\[5pt]
c(n)&=128(n-1)(n+1)^2.
\end{align*}
It is clear that $a(n)>0$ for any $n>3$.
By a routine computation, we get
\begin{align*}
b^2(n)-4a(n)c(n)=64(n^6+6n^5+15n^4+26n^3+25n^2+8n)>0, \quad \mbox{for $n>3$.}
\end{align*}
It suffices to show that
\begin{align}\label{eq-inequ}
\frac{V_n}{V_{n-1}}\geq \frac{-b(n)+\sqrt{b^2(n)-4a(n)c(n)}}{2a(n)}, \quad \mbox{for $n>3$.}
\end{align}
This inequality also implies the positivity of $V_n$ since its right-hand side is positive for any $n>3$.
(Note that $b(n)$ is negative.) However, it is difficult to directly prove \eqref{eq-inequ}. The key idea of our proof is to find an intermediate function $h(n)$ such that
\begin{align*}
\frac{V_n}{V_{n-1}}\geq h(n)\geq \frac{-b(n)+\sqrt{b^2(n)-4a(n)c(n)}}{2a(n)}, \quad \mbox{for $n>3$.}
\end{align*}
Let
\begin{align}\label{eq-h}
h(n)=\frac{16(n^3-n^2+1)}{n^3-n^2}, \quad \mbox{for $n\geq 2$,}
\end{align}
and we shall show that this function fulfills our purpose. This will be done in two steps.

First, we need to prove that
\begin{align}
h(n)-\frac{-b(n)+\sqrt{b^2(n)-4a(n)c(n)}}{2a(n)}\geq 0, \quad \mbox{for $n>3$.} \label{eq-inequ-1}
\end{align}
A straightforward computation shows that the quantity on the left-hand side is equal to
{\small{\begin{align*}
\frac{32(4n^6+7n^5+n^4+n^3+9n^2+8n+2)}
 {(n^4-n^2)(n+1)(n^5+2n^4+n^2+8n+4+(n^2-n)\sqrt{n^6+6n^5+15n^4+26n^3+25n^2+8n})},
\end{align*}}}
which is clearly positive for $n>3$.

Second, we need to prove that
\begin{align}\label{eq-inequ-2}
\frac{V_n}{V_{n-1}}\geq h(n), \quad \mbox{for $n>3$.}
\end{align} For the sake of
convenience, let $g(n)=\frac{V_n}{V_{n-1}}$. We use induction on $n$ to prove that
$g(n)\geq h(n)$ for $n>3$. By the recurrence relation \eqref{eq-rec}, we have
\begin{align}\label{eq-gr}
g(n+1)=\frac{8(3n^2+5n+1)}{(n+1)^2}-\frac{128(n-1)}{ng(n)},
\quad n\geq 1,
\end{align}
with the initial value $g(1)=8$.
It is clear that $g(3)=152/9=h(3)$ and $g(4)=625/38>49/3=h(4)$ by \eqref{eq-h} and \eqref{eq-gr}.
Assume that $g(n)>h(n)$, and we proceed to show that $g(n+1)>h(n+1)$.
Note that
\begin{align*}
 g(n+1)-h(n+1)
&=\frac{8(3n^2+5n+1)}{(n+1)^2}-\frac{128(n-1)}{ng(n)}
  -\frac{16(n^3+2n^2+n+1)}{n(n+1)^2}\\[5pt]
&=\frac{8(n^3+n^2-n-2)}{n(n+1)^2}-\frac{128(n-1)}{ng(n)}\\[5pt]
&=\frac{8(n^3+n^2-n-2)g(n)-128(n-1)(n+1)^2}{n(n+1)^2g(n)}.
\end{align*}
By the induction hypothesis, we have $g(n)>h(n)>0$ and thus
\begin{align*}
 g(n+1)-h(n+1) &>\frac{8(n^3+n^2-n-2)h(n)-128(n-1)(n+1)^2}{n(n+1)^2g(n)}\\[5pt]
&=\frac{128(2n^2-n-2)}{n^3(n-1)(n+1)^2g(n)}>0.
\end{align*}
Combining \eqref{eq-inequ-1} and \eqref{eq-inequ-2}, we obtained the inequality \eqref{eq-inequ}.
This completes the proof.
\qed

Wang and Zhu \cite[Theorem 2.1]{WangZhu} showed that if $\{z_n\}_{n\geq 0}$ is a log-concave sequence of positive integers with $z_0>1$, then $\{\sqrt[n]{z_n}\}_{n\geq 1}$ is strictly decreasing.
Applying their criterion to the Fennessey-Larcombe-French sequence, we obtain immediately the following result.

\begin{prop}
The sequence $\{\sqrt[n]{V_{n+1}}\}_{n\geq 1}$ is strictly decreasing.
\end{prop}

{\noindent \it Proof. } Let $\{z_n\}_{n\geq 0}$ be the sequence given by $z_n=V_{n+1}$.
It is clear that $z_0=V_1=8>1$. Moreover, by Theorem \ref{thm-main}, the sequence $\{z_n\}_{n\geq 0}$ is log-concave. Thus, $\{\sqrt[n]{z_n}\}_{n\geq 1}$ is strictly decreasing by \cite[Theorem 2.1]{WangZhu}.
\qed

\noindent{\bf Acknowledgements.} This work was supported by the 973 Project, the PCSIRT Project of the Ministry of Education and the National Science Foundation of China.

\end{document}